\input amstex
\documentstyle{amsppt}
\magnification=\magstep1
\hsize=5in
\vsize=7.3in
\TagsOnRight  
\topmatter
\title Algebraic Coset Conformal Field Theories II 
\endtitle
\author Feng Xu \endauthor

\address{Department of Mathematics, University of Oklahoma, 601 Elm Ave,
Room 423, Norman, OK 73019}
\endaddress
\email{xufeng\@ math.ou.edu}
\endemail
\abstract{   
Some mathematical questions relating to Coset  Conformal field theories 
(CFT) are considered in the framework of Algebraic Quantum Field Theory
as developed previously by us.  
 We  consider
the issue of fix point resolution in the diagonal coset of type A. We
show how to decompose certain reducible representations into irreducibles,
and  prove that the coset CFT gives rise to a unitary Modular Tensor Category
in the sense of Turaev, and therefore may be used to construct 3-manifold
invariants.
We prove that if the coset inclusion satisfies certain conditions which
can be checked in examples,  
the Kac-Wakimoto Hypothesis (KWH) is equivalent to the
the Kac-Wakimoto Conjecture (KWC), a result which seems to be
hard to prove by purely representation considerations. Examples are also
presented.}  
\endabstract        
\thanks 
I'd like to thank Professors J. Fuchs,   Karl-Henning 
Rehren and C. Dong for useful correspondences.    
1991 Mathematics Subject Classification. 46S99, 81R10. 
\endthanks  
            
\endtopmatter
\document   
            
\heading \S1. Introduction \endheading
This paper is a sequel to [X4]. Let us first recall some definitions
from [X4]. \par
Let $G$ be a simply connected  compact Lie group
and let $H\subset G$ be a Lie subgroup. Let $\pi^i$ be an irreducible
representations of $LG$ with positive energy at level $k$
\footnotemark\footnotetext{ When G is the direct product of simple groups,
$k$ is a multi-index, i.e., $k=(k_1,...,k_n)$, where $k_i\in \Bbb {N}$
corresponding to the level of the $i$-th simple group. The level of $LH$
is determined by the Dynkin indices of $H\subset G$. To save some writing
we write the coset as $H\subset G_k$.} 
on Hilbert space $H^i$.
Suppose when restricting to $LH$, $H^i$ decomposes as:
$$
H^i = \sum_\alpha H_{i,\alpha} \otimes H_\alpha 
$$, and  $\pi_\alpha$ are irreducible representations of $LH$ on
Hilbert space $H_\alpha$.  The set of $(i,\alpha)$ which appears in
the above decompositions will be denoted by $exp$. \par
We
shall use $\pi^1$ (resp. $\pi_1$) \footnotemark\footnotetext{ This
is slightly different from the notation $\pi^0$ (resp. $\pi_0$) in [X4]:
it seems
to be more appropriate since these representations correspond to identity
sectors.} to
denote the vacuum
representation of
$LG$ (resp. $LH$).
Let $A$ be the Vacuum Sector of the coset $G/H$ as defined on Page 5 of
[X4].  The decompositions above naturally give rise to a class of
covariant representations of $A$, denoted by $\pi_{i,\alpha}$ or simply
$(i,\alpha)$, by Th. 2.3 of [X4], $\pi_{1,1}$ is the Vacuum representation
of $A$. \par
In \S2.2 we consider the decompositions of certain reducible
representations in the diagonal cosets of type $A_{N-1}$ as considered in 
\S4.3 of [X4] when action of the Dynkin diagram automorphisms is 
not free (cf. Page 31 of [X4]),
which is part of the {\bf fixed  point resolution} 
problems known in physics literature (cf. [Gep], [LVW] and [SY]).
Such problems have been known for some time, and there are no even clear 
mathematical formulations of such questions before.  We will provide
further evidence that 
the results of [X4] provide the right conceptual framework for
understanding such questions. \par
            
We first prove a general lemma 2.1 which we 
believe will play an important role in all fixed point resolution
problems.  Using lemma 2.1 and lemma 2.2, we prove
(cf. (1) of Th. 2.3)  that certain $S$-matrices are non-degenerate.
It follows from Th.2.3 (Cor. 2.4) that 
the diagonal cosets of type $A_{N-1}$ 
give rise to a unitary Modular Tensor Category in the 
sense of Turaev (cf. P.74 and P.113 of [Tu]), and may be used 
to construct 3-manifold invariants (cf. P. 160 of [Tu]). We also 
calculate $S$ matrices when $N$ is prime.  The result agrees with 
some of the results of [FSS1], [SY] from different considerations. \par
To describe the results in \S3, let us  
denote by $S_{ij}$ (resp. $\dot{S_{\alpha \beta}}$) the $S$ matrices of
$LG$
(resp. $LH$) at level $k$ (resp. certain level of $LH$ determined by
the inclusion $H\subset G_k$). Define  \footnotemark\footnotetext {Our
$(j,\beta)$ corresponds to $(M,\mu)$ on P.186 of [KW], and it follows from
the definitions that $\langle (j,\beta), (1,1) \rangle$ is then equal
to $mult_M(\mu,p)$ which appears in 2.5.4 of [KW]. It is then easy to see
that our formula (1) is identical to 2.5.4 of [KW].}
$$
b(i,\alpha) = \sum_{(j,\beta) } {S_{ij}} \overline{\dot S_{\alpha \beta}}
\langle (j,\beta), (1,1) \rangle \tag 1
$$
Note the above summation is effectively over those $(j,\beta)$ such that 
$(j,\beta) \in exp$.
The Kac-Wakimoto Conjecture (KWC)  states that 
if $(i,\alpha) \in exp$, then $b(i,\alpha) >0$. \par
The Kac-Wakimoto Hypothesis (KWH) states that
if $\langle (j,\beta), (1,1) \rangle >0$ and $(i,\alpha) \in exp$,
then $S_{ij} \overline {\dot S_{\alpha \beta}} \geq 0$. \par
Note that since $S_{i1} >0$, $\dot S_{\alpha 1} >0, (1,1) \in exp$, 
KWH implies KWC. \par
KWC has proved to be true in all known examples. In
fact, in \S2.4 of [X4] an even stronger conjecture, Conjecture 2 (C2) is
formulated. \par 
Unfortunately KWH is not true. In [X2] counter examples were found
by using subfactors associated with conformal inclusions. However,  
KWH has checked to be true in so many examples, and it seems that it should
be true or 
equivalent to KWC 
under some general conditions. The first  main result in \S3.1 is 
 to describe such a condition (cf. Th.3.3). 
The condition \footnotemark\footnotetext{ If we identify $(1,\beta)$ with
$(M, \mu)$, where $M$ is the vacuum representation, as on P.186 of [KW],
then condition (2) is the statement that 
if $(M,\mu) \in S_m$, with $S_m$ defined on P.186 of [KW],
then $\mu$ must be the 
vacuum representation of the subalgebra.} is that :
$$
\text {\rm if}  \ \langle (1,\beta), (1,1) \rangle >0, \text {\rm then} 
 \ \beta =1 \tag 2
$$. 
Th. 3.3 states that if 
 $H\subset G_k$ satisfies (2), and 
 certain assumptions   in \S3.1 which are expected to be
 true in general ,
then KWH is equivalent to
KWC for the inclusion $H\subset G_k$. \par
Condition (2) can be shown to be equivalent to the normality of certain
inclusions, but we will not discuss this in this  paper.\par
In \S3.1 we also give an example which does not satisfy (2), and
verifies KWC but not KWH.  This is also the first example of non-conformal
inclusion which does not verify KWH. \par
It is interesting to note that Th. 3.3 
can be thought as  a statement about  representations of affine
Kac-Moody algebras without even mentioning von Neumann algebras, yet
it seems to be hard to obtain such results without using subfactor
theory (cf.[J]). We give another example of this nature in Prop.3.2. For more
such statements, see inequality on Page 11 of [X2] and in 
particular  (2) of Th.4.3 
of [X4].\par 
In \S3.2 we prove a property (Prop.3.4) of Conjecture 2 (C2) in [X4].
It states that if $H_1\subset H_2, H_2\subset G$ verifies C2, then 
$H_1\subset G$ also
verifies C2, thus reducing C2 to maximal inclusions which are classified
in [Dyn1], [Dyn2]. We also give an example related to $N=2$ superconformal
theories. \par 
\heading 2. Fixed point resolutions in the diagonal cosets of
type  $A_{N-1}$ \endheading 
\subheading {2.1 Preliminaries }
Let us first recall some definitions from [X2].
Let $M$ be a properly infinite factor
and  $\text{\rm End}(M)$ the semigroup of 
 unit preserving endomorphisms of $M$.  In this paper $M$ will always
be the unique hyperfinite $III_1$ factors. 
Let $\text{\rm Sect}(M)$ denote the quotient of $\text{\rm End}(M)$ modulo 
unitary equivalence in $M$. We  denote by $[\rho]$ the image of
$\rho \in \text{\rm End}(M)$ in  $\text{\rm Sect}(M)$.\par
 It follows from
\cite{L3} and \cite{L4} that $\text{\rm Sect}(M)$, with $M$ a properly
infinite  von Neumann algebra, is endowed
with a natural involution $\theta \rightarrow \overline \theta $  ;  
moreover,  $\text{\rm Sect}(M)$ is
 a semiring with identity denoted by $id$. \par
If given a normal
faithful conditional expectation
$\epsilon:
M\rightarrow \rho(M)$,  we define a number $d_\epsilon$ (possibly
$\infty$) by:
$$
d_\epsilon^{-2} :=\text{\rm Max} \{ \lambda \in [0, +\infty)|
\epsilon (m_+) \geq \lambda m_+, \forall m_+ \in M_+
\}$$ (cf. [PP]).\par
 We define
$$
d = \text{\rm Min}_\epsilon \{ d_\epsilon |  d_\epsilon < \infty \}.
$$   $d$ is called the statistical dimension of  $\rho$. It is clear
from the definition that  the statistical dimension  of  $\rho$ depends
only
on the unitary equivalence classes  of  $\rho$. 
The properties of the statistical dimension can be found in
[L1], [L3] and  [L4]. We will denote the statistical dimension of 
$\rho$ by $d_\rho$ in the following.  $d_\rho^2$ is called the minimal
index of $\rho$.\par
Recall from [X2] that  we denote 
by $\text{\rm Sect}_0(M)$ those elements
of 
$\text{\rm Sect}(M)$ with finite statistical dimensions.
For $\lambda $, $\mu \in \text{\rm Sect}_0(M)$, let
$\text{\rm Hom}(\lambda , \mu )$ denote the space of intertwiners from 
$\lambda $ to $\mu $, i.e. $a\in \text{\rm Hom}(\lambda , \mu )$ iff
$a \lambda (x) = \mu (x) a $ for any $x \in M$.
$\text{\rm Hom}(\lambda , \mu )$  is a finite dimensional vector 
space and we use $\langle  \lambda , \mu \rangle$ to denote
the dimension of this space.  $\langle  \lambda , \mu \rangle$
depends      
only on $[\lambda ]$ and $[\mu ]$. Moreover we have 
$\langle \nu \lambda , \mu \rangle = 
\langle \lambda , \bar \nu \mu \rangle $, 
$\langle \nu \lambda , \mu \rangle 
= \langle \nu , \mu \bar \lambda \rangle $ which follows from Frobenius 
duality (See \cite{L2} ).  We will also use the following 
notation: if $\mu $ is a subsector of $\lambda $, we will write as
$\mu \prec \lambda $  or $\lambda \succ \mu $.  A sector
is said to be irreducible if it has only one subsector. \par
Recall (cf.[L2]) for each $\rho\in \text{\rm End}$ and its conjugate
$\bar\rho$ with finite minimal index, there exists $R_\rho \in 
\text{\rm Hom}(id , \bar \rho \rho )$ and $ \bar R_\rho \in
\text{\rm Hom}(id , \rho \bar\rho )$ such that
$$
\bar R_\rho^* \rho(R_\rho) = 1,  R_\rho^* \bar \rho(\bar R_\rho) = 1
$$ and $||R_\rho||= ||\bar R_\rho||=\sqrt {d_\rho}$.  The minimal
left inverse $\phi_\rho$ of $\rho$ is defined by
$$
\phi_\rho (m) = R_\rho^* \bar\rho(m) R_\rho
$$.\par   
The following lemma plays a fundamental role in \S2.2.
\proclaim {Lemma 2.1}
Let $a,b,c\in End(M)$, $[c]=[ab]$ and
$a,b$ have finite statistical dimensions.  Suppose $\tau \in End(M)$ has
order
$t$ in $Sect(M)$ , i.e., $t$ is the least positive integer such that
$[\tau^t]=[id]$, and
$[a\tau] = [a], [\tau b]= [b]$. If 
$$
\langle c,c \rangle =t
$$, then Hom$(c,c)$ is an abelian algebra with dimension $t$ and hence
there exist irreducible sectors $c_1,...,c_t$ such that
$$
c=\sum_{1\leq k\leq t} c_k
$$.  Moreover, $d_{c_k}= \frac{1}{t} d_c , k=1,...,t.$
\endproclaim
{\it Proof:} 
From $[a\tau]=[a]$ we conclude by using Frobenius duality that
$\langle \bar a a, \tau^i \rangle \geq 1$, and since $\tau$ has order
$t$ in $Sect(M)$, we must have 
$$
\bar a a \succ \sum_{0\leq i \leq t-1} \tau^i
$$. Similarly 
$$
b \bar b  \succ \sum_{0\leq i \leq t-1} \tau^i
$$.
 
Since 
$$
\align           
t =\langle c,c\rangle  & = \langle ab,ab\rangle =
\langle \bar a a, b\bar b\rangle \\
&\geq \sum_{0\leq i\leq t-1}
\langle \bar a a, \tau^i \rangle \langle \tau^i, b\bar
b\rangle \geq t
\endalign
$$, it follows that all the $\geq$ are $=$ and in particular
$\langle \bar a a, id\rangle = \langle \bar b b, id\rangle =1$, i.e.,
both $a$ and $b$ are irreducible. 
It is enough to prove the case $c=ab$.
Since $ [\tau b]= [b]$, there exist unitary elements $v\in M$ such
that         
$$ v\in  Hom (b, \tau b)$$. 
Define $\tau_v = v^* \tau v$, then 
$$           
b= \tau_v b  
$$, and so   
$$           
b= \tau_v^t b 
$$. Since $[\tau^t]=[id]$, there exists a unitary $v_1$ such that
$$           
\tau_v^t = Ad_{v_1}
$$ and so $v_1 b = b v_1$.  Since $b$ is irreducible, $v_1$ is equal to
identity up to a complex number whose absolute value is 1, so
$$           
\tau_v^t =id 
$$.  From $[a\tau]=[a\tau_v]=[a]$ there exists a unitary $u$ such that
$ u\in  Hom (a \tau_v, a)$. It follows that 
$$           
u^t a =u^t a\tau_v^t = a u^t
$$, so $u^t = x.1$ with $x\in {\Bbb C}, $ $x$ has absolute value one, 
since $a$ is irreducible.
Define  $w=x^{\frac{1}{t}} u$ so that $w^t= 1$. Note $w\in Hom(ab,ab)$.
Denote by $\phi_a, \phi_b$ the minimal left inverse of $a,b$.
We claim that $ \phi_b \phi_a (w^i) = d_a d_b \delta_{i0}$ for $0\leq
i\leq t-1$.
It is enough to show 
$ \phi_b \phi_a (w^i) = 0$ for $0< i\leq t-1$.  Since $b$ is irreducible,
$\phi_a (w^i) \in \text {\rm Hom (b,b)} \equiv \Bbb {C}.1$, and so
$ \phi_b \phi_a (w^i) = R_{\bar b}^* \phi_a (w^i) R_{\bar b}$. But
$$
\align
R_{\bar b}^* \phi_a (w^i) R_{\bar b} &= 
R_{\bar b}^* R_a^*  \bar a(w^i) R_a  R_{\bar b} = R_a^* \bar a a(R_{\bar
b}^*) \bar a(w^i) R_a  R_{\bar b} \\
&= R_a^* \bar a ( a(R_{\bar  b}^*) w^i) R_a  R_{\bar b} = 
R_a^* \bar a (w^i a\tau_v^i(R_{\bar  b}^*)) R_a  R_{\bar b} \\
&= R_a^* \bar a (w^i) \bar a a (\tau_v^i(R_{\bar  b}^*))R_a  R_{\bar b} \\
&= R_a^* \bar a (w^i) R_a \tau_v^i(R_{\bar  b}^*) R_{\bar b}
\endalign
$$, and use the fact $\tau_v b =b$ it is easy to see that
$\tau_v^i(R_{\bar  b}^*) R_{\bar b} \in Hom (id, \tau_v^i)$.
Since  $ \langle \tau^i, id\rangle = 0, 0<i\leq t-1$, it follows that
$\tau_v^i(R_{\bar  b}^*) R_{\bar b}=0$ and so
 $ \phi_b \phi_a (w^i) = d_a d_b \delta_{i0}$ for $0\leq i\leq t-1$.
. \par 
If $\sum_{0\leq i\leq t-1} x_i w^i = 0, x_i\in {\Bbb C}$, 
multiplied both sides by
$w^{t-i}$ and apply $\phi_b \phi_a$, we get $x_i=0$. 
so $id, w, w^2, ..., w^{t-1}$ are linearily independent in $\text 
{\rm Hom(c,c)}$. and
since $\text {\rm Hom(c,c)}$ has dimension $t$, $\text {\rm Hom(c,c)}$ is
therefore an abelian
algebra with basis
$id, w, w^2, ..., w^{t-1}$.  Note $w^t=1$, so the minimal projections 
$P_k, k=1,...,t $ in $\text {\rm Hom(c,c)}$ are given by
$P_k= \sum_{0\leq j\leq t-1} \frac{1}{t} (\exp (\frac{2\pi ik}{t}) w)^j$.
Let $c_k\prec c$
be the irreducible sector corresponding to $P_k$, then by [L4]
$$           
d_{c_k} =  \phi_a\phi_b (P_k) = \frac{1}{t} d_ad_b
=\frac{1}{t} d_c, k=1,...,t.
$$.          
\hfill Q.E.D. \par
Next we will recall some of the results of [Reh] (also cf. [FRS]) and
introduce
notations. \par
Let $[\rho_i]$ denote the equivalence classes of irreducible
superselection
sectors \par
(cf.[GL])  in a finite set.  Suppose this set is closed under
conjugation 
and composition. We will denote the conjugate of $[\rho_i]$ by
$[\rho_{\bar i}]$
and identity sector by $[1]$ if no confusion arises, and let 
$N_{ij}^k = \langle [\rho_i][\rho_j], [\rho_k]\rangle $. We will
denote by $\{T_e\}$ a basis of isometries in $\text {\rm
Hom}(\rho_k,\rho_i\rho_j)$. 
The univalence of $\rho_i$ (cf. P.12 of [GL]) will be denoted by
$\omega_{\rho_i}$. \par
Let $\phi_i$ be the unique minimal 
left inverse of $\rho_i$, define:
$$ 
Y_{ij}:= d_{\rho_i}  d_{\rho_j} \phi_j (\epsilon (\rho_j, \rho_i)^*
\epsilon (\rho_i, \rho_j)^*), \tag 0
$$ where $\epsilon (\rho_j, \rho_i)$ is the unitary braiding operator
 (cf. [GL] ). \par
We list two properties of $Y_{ij}$ (cf. (5.13), (5.14) of [Reh]) which
will be used in \S2.2:
$$
\align
Y_{ij} = Y_{ji} & = Y_{i\bar j}^* = Y_{\bar i \bar j} \tag 1 \\
Y_{ij} = \sum_k N_{ij}^k \frac{\omega_i\omega_j}{\omega_k} d_{\rho_k} \tag
2
\endalign
$$.  Let us explain the proof of (2) since similar but different proof
appears in the proof of Lemma 2.2. \par
We have:
$$
\align
\phi_j (\sum_e T_e T_e^*\epsilon (\rho_j, \rho_i)^*
\epsilon (\rho_i, \rho_j)^*  T_e T_e^*) & = 
\sum_e \frac{\omega_i\omega_j}{\omega_k} \phi_j (T_e T_e^*) \\
&= \sum_e \frac{\omega_i\omega_j}{\omega_k} N_{ij}^k
\frac{d_{\rho_k}}{d_{\rho_i}d_{\rho_j}}
\endalign
$$, where in the first $=$ we used the Monodromy equation (cf. [FRS]
or [X1]), and the second $=$ follows from [L2].\par
 Define 
$\tilde \sigma := \sum_i d_{\rho_i}^2 \omega_{\rho_i}^{-1}$.
If the matrix $(Y_{ij})$ is invertible,
by Proposition on P.351 of [Reh] $\tilde \sigma$ satisfies
$|\tilde \sigma|^2 = \sum_i d_{\rho_i}^2$. 
Suppose $\tilde \sigma= |\tilde \sigma| \exp(i x), x\in {\Bbb R}$.
Define matrices 
$$
S:= |\tilde \sigma|^{-1} Y, T:=  \exp(i \frac{x}{3}) Diag(\omega_{\rho_i})
\tag 3
$$.  Then these matrices satisfy the algebra:
$$
\align
SS^{\dag} & = TT^{\dag} =id, \tag 4  \\
TSTST&= S, \tag 5 \\
S^2 =C, TC=CT=T, \tag 6
\endalign
$$ 
where $C_{ij} = \delta_{i\bar j}$ is the conjugation matrix. Moreover
$$
N_{ij}^k = \sum_m \frac{S_{im} S_{jm} S_{km}^*}{S_{1m}}. \tag 7
$$
(7) is known as Verlinde formula. \par
Now let us consider an example which verifies (1) to (7) above.
Let $G= SU(N)$. We denote $LG$ the group of smooth maps
$f: S^1 \mapsto G$ under pointwise multiplication. The
diffeomorphism group of the circle $\text{\rm Diff} S^1 $ is 
naturally a subgroup of $\text{\rm Aut}(LG)$ with the action given by 
reparametrization. In particular the group of rotations
$\text{\rm Rot}S^1 \simeq U(1)$ acts on $LG$. We will be interested 
in the projective unitary representation $\pi : LG \rightarrow U(H)$ that 
are both irreducible and have positive energy. This means that $\pi $ 
should extend to $LG\ltimes \text{\rm Rot}\ S^1$ so that
$H=\oplus _{n\geq 0} H(n)$, where the $H(n)$ are the eigenspace
for the action of $\text{\rm Rot}S^1$, i.e.,
$r_\theta \xi = \exp^{i n \theta}$ for $\theta \in H(n)$ and 
$\text{\rm dim}\ H(n) < \infty $ with $H(0) \neq 0$. It follows from 
\cite{PS} that for fixed level $K$ which
is a positive integer, there are only finite number of such 
irreducible representations indexed by the finite set
$$           
 P_{++}^{h}  
= \bigg \{ \lambda \in P \mid \lambda 
= \sum _{i=1, \cdots , N-1}
\lambda _i \Lambda _i , \lambda _i \geq 1\, ,
\sum _{i=1, \cdots , n-1}
\lambda _i < h \bigg \}
$$           
where $P$ is the weight lattice of $SU(N)$ and $\Lambda _i$ are the 
fundamental weights and $h=N+K$.  We will use 
$1$ to denote the trivial representation  of 
$SU(N)$. For $\lambda , \mu , \nu \in  P_{++}^{K}$, define
$$
N_{\lambda \mu}^\nu  = \sum _{\delta \in P_{++}^{K} } \frac{S_{\lambda
\delta} 
S_{\mu \delta} S_{\nu \delta}^*}{S_{1\delta}} \tag 8
$$ 
where $S_{\lambda\delta}$ is given 
by the Kac-Peterson formula:
$$           
S_{\lambda \delta} = c \sum _{w\in S_N} \varepsilon _w \exp
(iw(\delta) \cdot \lambda 2 \pi /n) \tag 9
$$.           
Here  $\varepsilon _w = \text{\rm det}(w)$ and $c$ is a normalization 
constant fixed by the requirement that $(S_{\lambda\delta})$ 
is an orthonormal system. 
It is shown in \cite{Kac} P.288 that $N_{\lambda \mu}^\nu $ are
non-negative 
integers. Moreover, define $ Gr(C_K)$ 
to be the ring whose basis are elements 
of $ P_{++}^{K}$ with structure constants $N_{\lambda \mu}^\nu $.
  The natural involution $*$ on $ P_{++}^{K}$ is 
defined by $\lambda \mapsto \lambda ^* =$ the conjugate of $\lambda $ as 
representation of $SU(N)$.  All the irreducible representations of 
$Gr(C_K)$ are given by $\lambda \rightarrow
\frac{S_{\lambda\mu}}{S_{1\mu}}$ for some $\mu$. \par
The irreducible positive energy representations of $ L SU(N)$ at level
$K$ give rise to an irreducible conformal precosheaf ${\Cal A}$ 
and 
its covariant representations (cf. P.362 of [X1]). The unitary equivalent
classes of such  representations are the superselection sectors.
  We will use 
$\lambda$ to denote such   representations. \par
For $\lambda$ irreducible, the univalence $\omega_\lambda$ is given by
an explicit formula .  
Let us first 
define 
$$
\Delta_\lambda = 
\frac {c_2(\lambda)}{K+N} \tag 10
$$ where $c_2(\lambda)$ is the value of
Casimir
operator on representation of $SU(N)$ labeled by dominant weight
$\lambda$ (cf. 1.4.1 of [KW]).
 $\Delta_\lambda$ is usually called the conformal dimension.
Then
we have 
$\omega_\lambda = \exp({2\pi i} \Delta_\lambda)$. \par

Define the central charge (cf. 1.4.2 of [KW]) 
$$
C_G := \frac {K \text {\rm dim(G)}}{K+N} \tag 11
$$ and $T$ matrix as
$$
T=diag(\dot\omega_\lambda) \tag 12
$$, where $\dot\omega_\lambda = \omega_\lambda exp (\frac{-2\pi i
C_G}{24})
$.  By Th.13.8 of [Kac] $S$ matrix as defined in (9) and $T$ matrix
in (12) satisfy relation (4), (5) and (6). \par
By Cor.1 in \S34 of [W2],  The fusion ring generated by all
$\lambda \in   P_{++}^{(K)}$
is isomorphic to $ Gr(C_K)$, with structure constants $N_{\lambda
\mu}^\nu$ as defined in (8).  One may therefore ask what are the $Y$ matrix
(cf. (0)) in this case. By using (2) and the formula for 
 $N_{\lambda
\mu}^\nu$, a simple calculation shows:
$$
Y_{\lambda \mu} = \frac{S_{\lambda \mu}}{S_{1\mu}}
$$, and it follows that $Y_{\lambda \mu}$ is nondegenerate, and $S,T$
matrices
as defined in (3) are indeed the same $S,T$ matrix defined in (8) and
(11),
which is a surprising fact.  If the analogue of Cor.1 in \S34 of 
[W2] is established for other types of simple and simply connected Lie 
groups, then this fact is also
true for other types of groups by the same argument. \par
In \S2.2 we will also consider the case when $G$ is the direct product
of two type $A$ groups. In that case the $S,T$ matrices are just the
tensor product of  the $S,T$ matrices corresponding to each subgroup.
\subheading {2.2 Fixed point resolutions}
We preserve the set up of \S4.3 of [X4].  
We consider the coset 
$G:=SU(N)_{m'} \times SU(N)_{m''}/H:=SU(N)_{m'+m''}$, where
the embedding $H\subset G$ is diagonal. 
Let $\Lambda_1,...,\Lambda_{N-1}$ be the fundamental weights of $SL(N)$.
Let $k\in \Bbb N$. Recall that the set of integrable weights of the affine
algebra $\widehat{SL(N)}$ at level $k$ is the following subset of the
weight
lattice of  $SL(N)$:
$$
P_{++}^{(h)} = \{ \lambda = \lambda_1 \Lambda_1 +...+ \lambda_{N-1}
 \Lambda_{N-1} | \lambda_i\in \Bbb N,  \lambda_1+...+ \lambda_{N-1} <h
\}
$$ where $h=k+N$.
 This set admits a $\Bbb {Z}_N$ automorphism generated  
by
$$
\sigma_1: \lambda=( \lambda_1, \lambda_2,..., \lambda_{N-1}) \rightarrow
\sigma_1( \lambda) = 
(h-\sum_{j=1}^{N-1} \lambda_j, \lambda_1,...,\lambda_{N-2})
$$.  We define the color $\tau(\lambda):\equiv \sum_i (\lambda_i-1)i 
 \text{\rm mod} (N)$ and $Q$ to be the root lattice of $\widehat{SL(N)}$
(cf. \S1.3 of [KW]). Note that $\lambda\in Q$ iff
$\frac{1}{N}\tau(\lambda) \in {\Bbb Z}$. 
\par 
We
use $i$ (resp.$\alpha$) to denote the irreducible positive
energy representations of $LG$ (resp.$LH$).  To
compare our notations with that \S2.7 of 
[KW], note that our $i$ is $(\Lambda',\Lambda'')$
of [KW] , and our $\alpha$ is $\Lambda$ of [KW]. 
We will  identify $i=(\Lambda',\Lambda'')$ and $\alpha =\Lambda$
where $\Lambda',\Lambda''$, $\Lambda$ are the weights of 
$SL(N)$ at levels $m', m'', m'+m''$ respectively.
Suppose 
$$
i=({\Lambda_1}',{\Lambda_1}''), j=({\Lambda_2}',{\Lambda_2}''),
k=({\Lambda_3}',{\Lambda_3}''), \alpha= {\Lambda_1}, \beta={\Lambda_2},
\delta=\Lambda_3.
$$           
Then the fusion coefficients
$N_{ij}^k := N_{{\Lambda_1}'{\Lambda_2}'}^{\Lambda_3'}
N_{{\Lambda_1}''{\Lambda_2}''}^{{\Lambda_3}''}$  
(resp. $N_{\alpha\beta}^\delta := N_{\Lambda_1\Lambda_2}^{\Lambda_3} $ 
)of $LG$ (resp. $LH$)
are given by Verlinde formula (cf. \S2.1).
Recall $\pi_{i,\alpha}$ are the covariant representations of
the coset $G/H$. The set of all $(i,\alpha):=
(\Lambda',\Lambda'',\Lambda)$ 
which appears
in the decompositions of $\pi^i$ of $LG$ with respect to $LH$
is denoted by $exp$. This set is determined on P.194 of
[KW] to be $(\Lambda',\Lambda'',\Lambda) \in exp$ iff 
$\Lambda'+\Lambda''-\Lambda \in Q$. 
The  ${\Bbb Z_N}$ action on $(i,\alpha), \forall i,
\forall \alpha$  is denoted by $\sigma(i,\alpha):= 
(\sigma (\Lambda'), \sigma (\Lambda''), \sigma (\Lambda))$,
$\sigma \in  {\Bbb Z_N}$. This is also known as diagram automorphisms
since they corresponds to the automorphisms of Dykin diagrams.  Note that
this  ${\Bbb Z_N}$ action preserves $exp$ and therefore induces a 
${\Bbb Z_N}$ action on  $exp$. 
\par
We define a vector space $W$ over $\Bbb {C}$ whose 
orthonormal basis are denoted by
$i\otimes \alpha$ with $i=(\Lambda',\Lambda'') , \alpha=\Lambda$. $W$
is also a commutative ring with structure constants given by
$N_{ij}^k N_{\alpha\beta}^\delta$. Let 
$V$ be the vector space over $\Bbb {C}$ whose basis are given by
the irreducible components of $\sigma_i a_{1\otimes \bar \alpha}$
(cf. \S4.3 of [X4]) . 
Then $V= V_0 \oplus V_1$, where $V_0$ is a subspace of $V$ whose 
 basis are given by
the irreducible components of $\sigma_i a_{1\otimes \bar \alpha}$ 
with $(i,\alpha)\in exp$, and $V_1$ is the orthogonal complement of
$V_0$ in $V$. The composition of sectors gives  $V$  a ring structure.
   By (1) of theorem 4.3 of [X4], the irreducible
subrepresentations of $(i,\alpha)$ of the coset are in one-to-one
correspondence with the basis of $V_0$ and this map is a ring
isomorphism by (1) of Prop.4.2 of [X4], and we will identify 
 the irreducible
subrepresentations of $(i,\alpha)$ of the coset with the basis of 
$V_0$ in the following when no confusion arises.  Note that $V_0$ is 
a subring of $V$ and $V_0.V_1\subset V_1$. \par
Define a linear map $P: W\rightarrow V$ such that $P(i\otimes \alpha) =
\sigma_i a_{1\otimes  \bar \alpha}$. By Th.4.3 of [X4] 
$$
P(i\otimes \alpha) =   
\sigma_i a_{1\otimes  \bar\alpha} = P(i'\otimes \bar\alpha') =
\sigma_{i'} a_{1\otimes  \alpha'}      
$$ iff $\sigma^s(i)=i', \sigma^s(\alpha)=\alpha'$ for some $s\in
{\Bbb Z}$.  Also 
$\langle P(i\otimes \alpha), P(j\otimes \beta) \rangle= 0
$ if  $P(i\otimes \alpha)\neq P(j\otimes \beta)$ by (*) of \S4.3 of
[X4].

Note $P(\sigma (1) \otimes \sigma (1)) = 1$ and 
$P$ is  a ring homomorphism from  $W$ to $V$. 
Define $W_0:= P^{-1} (V_0), W_1:= P^{-1} (V_1)$, then
$W= W_0 \oplus W_1$ since $exp$ is $\sigma$ invariant.  Note that
$i\otimes \alpha \in W_0$ iff $i- \alpha \in Q$.
Define the action of $Z_N$ on $W$ as $\sigma(i\otimes \alpha)
= \sigma(i) \otimes \sigma(\alpha)$. 
Much of the
following depends on the relation between $W$ and $V$. \par

Assume
$\sigma^s(i\otimes \alpha) = i\otimes \alpha$ for some $i\otimes \alpha
\in W$, 
and $0<s\leq N$ is the least positive integer with this property.  Let
$t=\frac{N}{s}$. By equation (*) on Page 30 of [X4] we have:
$$
\langle P(i\otimes \alpha), P(i\otimes \alpha) \rangle =t
$$.  Our first question is to decompose $P(i\otimes \alpha)$ when 
$t>1$. \par
Apply lemma 2.1 to the present case with 
$P(i\otimes 1) =a, P(1\otimes \alpha) =b , P(\sigma^s \otimes 1)=\tau$, we
conclude that
there exists
$c_1,...,c_t \in V$  such that
$$
P(i\otimes \alpha)= \sum_{1\leq k\leq t} c_k
$$ and $d_{c_k} = \frac{1}{t}d_i d_\alpha, k=1,...,t
$.  Note that if $P^{-1}(P(i\otimes \alpha))= \{i_1\otimes \alpha_1,...
i_s\otimes \alpha_s \}$, then $st = N$. \par
Note we identify the covariant representations of the coset with 
the basis of $P(W_0) = V_0$. 
The univalence of $A:=P(i\otimes \alpha), i\otimes \alpha \in W_0$ are
given
by:
$\omega_A = \exp(2\pi i (\Delta_i - \Delta_\alpha))$,
 where
$\Delta_i, \Delta_\alpha$ are the conformal dimensions (cf. \S2.1, and 
if
$i=(\Lambda',\Lambda''), \Delta_i:= \Delta_{\Lambda'} + \Delta_{\Lambda''}
$). 
Note
if $A\succ a$, then $ \omega_a = \omega_A$.   
The univalence is only defined for covariant sectors which 
correpond to elements of $V_0$.  However, for convenience 
let  us define  $\omega_{i\otimes  \alpha} := \omega_i
\omega_\alpha^{-1}$ for $i\otimes  \alpha \in W$.
Then if  $i\otimes \alpha \in
W_0$,  $\omega_i \omega_\alpha^{-1}$ 
which is the univalence of $P(i\otimes
\alpha)$ depends only on $P(i\otimes
\alpha)$ , i.e.,  $\omega_i \omega_\alpha^{-1}$ depends only on 
the orbit of $i\otimes \alpha$ under the $\Bbb {Z_N}$ action. \par
Suppose $A:=P(i\otimes \alpha) , B:=P(j\otimes \beta) , i\otimes \alpha, 
j\otimes \beta \in W_0$.
Let $\phi_A, \phi_B$ be
the
unique minimal left inverses of $A,B$.  Define
$$
Y_{AB}:= d_A d_B \phi_B \phi_A (\epsilon (B,A)^* \epsilon (A,B)^*) \tag 1
$$ . Note (1) is similar to (0) of \S2.1 , the difference here is that
our $A,B$ may be reducible and hence we need to include $\phi_B$
in the definition since $\phi_A (\epsilon (B,A)^* \epsilon (A,B)^*)$
may not be a scalar.   \par
To avoid confusions we will denote $S$ matrices associated to 
indices $i$
(Recall from \S2.1 this is the tensor product of $S$ matrices
associated to two type A subgroup of $G$) by $S_{ij}$ and the 
 $S$ matrices associated to indices $\alpha$ by 
$\dot{S_{\alpha \beta}}$.
\proclaim {Lemma 2.2}
Suppose $A:=P(i\otimes \alpha) , B:=P(j\otimes \beta) , i\otimes \alpha,
j\otimes \beta \in W_0$, and $A=\sum_{1\leq i \leq t} c_i$ with 
$d_{c_i} = \frac{1}{t} d_A$. Then: \par
(1) $\langle c_i B, P(k\otimes \delta) \rangle = \frac{1}{t}
 \langle A B, P(k\otimes \delta) \rangle, \forall k\otimes \delta 
\in W$ ; \par
(2) $Y_{AB} = \frac{S_{ij}}{S_{11}} \frac{\overline {\dot S_{ \alpha 
\beta}}}{\dot S_{11}} $;\par
(3) $Y_{c_iB} = \frac{1}{t}  Y_{AB}$;\par
(4) If $B= \sum_j b_j$, then $\sum_j Y_{c_i b_j} = Y_{c_iB}.$ 
\endproclaim
{\it Proof:}
Ad (1): Denote by $C:=  P(k\otimes \delta)$. Then by Frobenius duality
$$
\langle c_i B, C \rangle = \langle c_i,  C\bar B \rangle
$$. By the definitions of $B,C$, 
$$
 C\bar B = l A + \sum D
$$, where $l\geq 0$ is an integer, and $D$ are elements of the form
$P(k'\otimes \delta')$ which 
are  different from $A$, and by $(*)$ in \S4.3 of [X4] 
$\langle D, A \rangle = 0$, and so 
$\langle D, c_i \rangle = 0$.  It follows that
$$
\langle c_i B, C \rangle = \langle c_i, lA \rangle = l
$$, which is independent of $i$.  Since $A=\sum_{1\leq i\leq t} c_i$,
(1) follows. \par
Ad (2): The main point of the proof is that even though $A,B$ may be reducible,
their univalence  are complex numbers, so the monodromy
equation (cf. [FRS] or Page 359 of [X1]) holds,  and we  have:
$$
\align
 \phi_A \phi_B (\epsilon (B,A)^* \epsilon (A,B)^*) & = 
\phi_A \phi_B (\sum_{e \in V} T_e T_e^* \epsilon (B,A)^* \epsilon
(A,B)^*)T_e
T_e^*) \\
& = \sum_{e\in V} \frac{\omega_A \omega_B}{\omega_e}N_{AB}^e
\frac{d_e}{d_A
d_B}
\\
&=\sum_{e\in V} \frac{\omega_i \omega_j \omega_{\alpha}^{-1}
\omega_{\beta}^{-1}}{\omega_e} \langle AB,e\rangle  \frac{d_e}{d_A d_B}  
\endalign
$$, where $e\in V$ means we sum over the basis of $V$. Note that
the summation above is effectively over $V_0$ since $A\in V_0, B\in V_0$
and $V_0$ is a subring.   Suppose
$P(k\otimes \delta) =e_1+e_2+...+e_m, k\otimes \delta \in W_0$,
with
$d_{e_i} = \frac{1}{m} 
d_k d_\delta$ , $\omega_{e_i} = \omega_k \omega_\delta^{-1}$ . Assume
$P^{-1}(P(k\otimes  \delta)) = k_1 \otimes  \delta_1,
..., k_n \otimes  \delta_n$, with $mn=N$.  Then
$$
\sum_{e_i} \frac{\omega_i \omega_j \omega_{\alpha}^{-1}
\omega_{\beta}^{-1}}{\omega_{e_i}} \langle AB,e_i\rangle
\frac{d_{e_i}}{d_A
d_B}
= \frac{1}{m} \frac{\omega_i \omega_j \omega_{\alpha}^{-1}
\omega_{\beta}^{-1}}{\omega_k \omega_\delta^{-1}} \langle AB, P(k\otimes 
\delta)\rangle \frac{d_k d_\delta}{d_A
d_B}
$$. From
$$
\align
\langle AB, P(k\otimes
\delta)\rangle & =  \langle \sigma_i a_{1\otimes \bar \alpha} 
\sigma_j a_{1\otimes \bar \beta}
, \sigma_k a_{1\otimes \bar \delta}\rangle =
N_{ij}^{k'} N_{\bar \alpha \bar \beta}^{\bar \delta'} 
\langle \sigma_{k'} a_{1\otimes \bar \delta'}, \sigma_k a_{1\otimes \bar
\delta} \rangle \\
& = m (N_{ij}^{k_1} N_{\bar \alpha \bar \beta}^{\bar
\delta_1} +... + N_{ij}^{k_n} N_{\bar \alpha \bar \beta}^{\bar
\delta_n})
\endalign
$$, cf. equation (*) on Page 30 of [X4], we conclude that
$$
\sum_{e_i} \frac{\omega_i \omega_j \omega_{\alpha}^{-1}
\omega_{\beta}^{-1}}{\omega_e} \langle AB,e_i\rangle  \frac{d_{e_i}}{d_A
d_B}
=  \frac{\omega_i \omega_j \omega_{\alpha}^{-1}
\omega_{\beta}^{-1}}{\omega_k \omega_\delta^{-1}} 
(N_{ij}^{k_1} N_{\bar \alpha \bar \beta}^{\bar
\delta_1} +... + N_{ij}^{k_n} N_{\bar \alpha \bar \beta}^{\bar
\delta_n})
\frac{d_k d_\delta}{d_A
d_B}
$$, and so 
$$
\align
\phi_A \phi_B (\epsilon (B,A)^* \epsilon (A,B)^*) 
&=\sum_{k\otimes  \delta \in W} \frac{\omega_i \omega_j
\omega_{\bar\alpha}^{-1}
\omega_{\bar\beta}^{-1}}{\omega_{k} \omega_{\bar\delta}^{-1}}
 N_{ij}^{k}
N_{\bar \alpha\bar\beta}^{\bar\delta} \frac{d_k d_\delta}{d_i
d_\alpha d_j d_\beta}
\\
&=\sum_{k} \frac{\omega_i \omega_j}{\omega_{k}} N_{ij}^{k}
\frac{d_{k}}{d_id_j}
\times \sum_{\delta}
\frac{\omega_\alpha^{-1}\omega_{\beta}^{-1}}{ \omega_{\delta}^{-1}}
N_{\alpha\beta}^{\delta} \frac{ d_{\delta}}{ d_\alpha d_\beta}
\\
& = \frac{S_{ij}}{S_{11}} \frac{\overline {\dot S_{\alpha \beta}}}
{\dot S_{11}} 
\endalign   
$$, where in the last $=$ we have used (2) and the comment after (12)
in   \S2.1.\par
Ad (3): As in the proof of 
(2) let  $P(k\otimes \delta), k\otimes \delta
\in
W_0$,
$P(k\otimes \delta) =e_1+...+e_m, P^{-1}(P(k\otimes \delta))=
\{ k_1\otimes \delta_1,..., k_n\otimes \delta_n \}$,
with $mn=N$.  Then
$$
\align
\sum_{e_j} \langle c_iB, e_j\rangle \frac{\omega_{c_i} \omega_B}
{\omega_{e_j}} d_{e_j} & =  \frac{1}{m}\langle c_iB,
P(k\otimes
\delta)\rangle
\frac{\omega_{A}
\omega_B}
{\omega_{k\otimes \delta}} d_{k} d_\delta \\
&= \frac{1}{t} \frac{1}{m}\langle AB, P(k\otimes
\delta)\rangle
\frac{\omega_{A}
\omega_B}
{\omega_{k\otimes \delta}} d_{k} d_\delta \\ 
&= \frac{1}{t} \sum_{e_j} \langle AB, e_j\rangle \frac{\omega_{A}
\omega_B}
{\omega_{e_j}} d_{e_j} 
\endalign
$$, where on the second $=$ we used (1), and 
(3)  follows immediately.
\par
Ad (4): First note $\omega_{b_j} =\omega_B$. By (2) of \S2.1, we have:
$$
\align
\sum_j Y_{c_i b_j} & = 
\sum_j \sum_e \langle c_ib_j, e\rangle \frac{\omega_{c_i}
\omega_{b_j}}{\omega_e} d_e=  
\sum_e \sum_j \langle c_ib_j, e\rangle \frac{\omega_{c_i}
\omega_B}{\omega_e} d_e \\
&= \sum_e  \langle c_i B, e\rangle \frac{\omega_{c_i}
\omega_B}{\omega_e} d_e = Y_{c_i B}
\endalign
$$. \par

\hfill Q.E.D. \par 

Recall the basis of $V_0$ corresponds to a finite set of 
irreducible covariant sectors of the coset: it is closed under 
composition, conjugation and contains identity by Th.4.3 of 
[X4]. Define the $Y$ matrix as in (0) of \S2.1. Then we have:
\proclaim{Theorem 2.3}
(1).The matrix $Y$ is invertible. \par
(2) The number $\tilde \sigma:= \sum_{e\in V} d_e^2 \omega_e^{-1}$
is given by 
$$
\tilde \sigma = \frac{1}{N} S_{11}^{-1} \dot{{S_{11}}}^{-1}
\exp (\frac{-6\pi i}{24} (C_G -C_H))
$$, where 
$$
\frac{1}{24}(C_G - C_H) =
\frac{N^2-1}{24} - \frac{N(N^2-1)}{24} ( \frac{1}{m'+N} +
\frac{1}{m'+N} - \frac{1}{m'+m''+N})
$$.  
\endproclaim
{\it Proof:}
By (i) of Prop. on Page 351 of [Reh] it is enough to show that
if $e$ is such that
$$
Y_{eg} = d_{e}d_{g} \ \forall g 
$$, then $e=1$. Suppose $B:=P(j\otimes \beta)=\sum_{1\leq i\leq t_B} g_i
\succ g, B\in V_0$, by (4) of lemma 2.2
$$
Y_{eB} = \sum_i Y_{eg_i} = \sum_i d_e d_{g_i} = d_e d_B
$$. Suppose  $ A=P(i\otimes \alpha) = \sum_{1\leq i\leq t} c_i \succ e$, 
$d_{c_i} = 
\frac{1}{t} d_c.$  
Then by (3) of Lemma 2.2 we have:
$$
Y_{AB}= d_A d_B,
\forall B\in V_0
$$  and by (2) of lemma 2.2  we have
$$
\frac{S_{ij}}{S_{11}}\frac{\overline {\dot S_{\alpha \beta}}}{\dot S_{11}}
= d_i d_j d_\alpha d_\beta
$$, and it follows that
$$
\frac{S_{ij}}{S_{1j}} \frac{\overline {\dot S_{ \alpha\beta}}}{\dot S_{
1\beta}} = d_i
d_\alpha
$$ for any $j\otimes \beta\in W_0$. Note 
$$
|\frac{S_{ij}}{S_{1j}}| \leq d_i, |\frac{\overline {\dot S_{ \alpha\beta}}}
 {\dot S_{
1\beta}}| \leq d_\alpha
$$, there must exist $a_{ij} \in {\Bbb R}$ such that
$$
\frac{S_{ij}}{S_{1j}} = \exp(ia_{ij}) d_i,
\frac{ \dot S_{ \alpha\beta}}{\dot S_{
1\beta}} = \exp(ia_{ij}) d_\alpha
$$, for any $j\otimes \beta\in W_0$. Suppose
$(j,\beta)=(\Lambda',\Lambda'';
\Lambda)$, $(i,\alpha)=(
M',M''; M)$, then $\Lambda'+\Lambda''-\Lambda \in Q, M'+M''-M\in Q.$
From 
$$
\frac{ \dot S_{ \alpha\beta}}{\dot S_{
1\beta}} = \exp(ia_{ij}) d_\alpha
$$ we get
$$
|\frac{ S_{ M\Lambda}}{S_{
1\Lambda}}| =  d_M
$$ for any $\Lambda \in P_+^{(K'+K'')}$, and so
$$
\sum_\Lambda | S_{ M\Lambda}|^2 =  d_M^2 \sum_\Lambda | S_{
1\Lambda}|^2 
$$, i.e., $d_M^2=1$. It follows that $M\bar M$ is the Vacuum sector, and so
$M\Lambda$ is always irreducible.  Choose $\Lambda$ corresponding 
to the defining representation of $SU(N)$, it follows from the fusion
rules that $M$ must be of the form $\sigma^{s}(1)$ for some $s\in {\Bbb
Z}$. Since $P(i\otimes \alpha) = P(\sigma^{-s}(i \otimes \alpha)) =A$,  
replacing $(i,\alpha)= (M',M''; \sigma^{s}(1))$ by 
$\sigma^{-s} (i,\alpha)= (\sigma^{-s}(M'),\sigma^{-s}( M'');1)
$ if necessary, we may assume $s=0$, i.e.,  $\alpha =1$. \par
Similarly $(M',M'')= (\sigma^{s_1}(1), \sigma^{s_2}(1))$, and using 
$\alpha =1$ we have
$$
\frac{ S_{ M'\Lambda'}}{S_{
1\Lambda'}} \frac{ S_{ M''\Lambda''}}{S_{
1\Lambda'}}=  d_{M'}d_{M''}, \forall (\Lambda', \Lambda'')
$$. Choose $\Lambda'$ to be the Vaccum representation and $ \Lambda''$
corresponding to  the defining representation of $SU(N)$ and use (2.2.15)
of [KW], we conclude $M'$ is the Vacuum representation, and similarly
$M''$ is also the Vacuum representation.\par
So we have proved $P(i\otimes \alpha)\succ e$ is the Vacuum sector, and
therefore
$e$
must be the Vacuum sector.\par
Ad (2): First we claim:
$$
\sum_{e\in V_0} d_e^2 \omega_e^{-1}
= \frac{1}{N} \sum_{g\in W_0} d_g^2 \omega_g^{-1}
$$. Suppose $P(g)=f$, $P^{-1}(P(g)) = \{g_1,...,g_n \}$, and
$P(g) = f=e_1+...+e_m$, with $mn=N$.  Since $\omega_{e_i} = \omega_g$,
$d_{e_i} = \frac{d_{g}}{m}, d_{g_j}= d_g, \omega_{g_j} = \omega_g$, we have 
$$
\sum_{e_i} d_{e_i}^2 \omega_{e_i}^{-1} 
= \sum_{e_i} \frac{d_{g}^2}{m^2}  \omega_{g}^{-1} 
= \frac{d_{g}^2}{m}  \omega_{g}^{-1}
$$, and so 
$$
\sum_{e_i} d_{e_i}^2 \omega_{e_i}^{-1}
= \frac{1}{N} n {d_{g}^2}  \omega_{g}^{-1}
= \frac{1}{N} \sum_{g_j} {d_{g_j}^2}  \omega_{g_j}^{-1}
$$, and it follows that
$$
\sum_{e\in V_0} d_e^2 \omega_e^{-1}
= \frac{1}{N} \sum_{g\in W_0} d_g^2 \omega_g^{-1}
$$.  Next let us show
$$
\sum_{g\in W_1} d_g^2 \omega_g^{-1} = 0
$$.  Again Suppose $ g=i\otimes \alpha, P(g)=f$, $P^{-1}(P(g)) =
\{g_1,...,g_n \}$, and
$P(g) = f=e_1+...+e_m$, with $mn=N$.  So we have 
$\sigma^n (i,\alpha) = (i,\alpha), g_k = \sigma^{k-1} g_1,1\leq k\leq n$.
Note that
$$
\omega_{\sigma (i)} \omega_{\sigma (\alpha)}^{-1} = \exp (\frac{-2\pi
i (\tau(i)- \tau(\alpha))}{N}) \omega_i \omega_\alpha^{-1}
$$. Denote by 
$z:= \exp (\frac{2\pi
i (\tau(i)- \tau(\alpha))}{N})$, then $z^n=1$ since 
$\sigma^n (i,\alpha) = (i,\alpha)$, but $z\neq 1$ since $i\otimes \alpha$
is not in $W_0$, i.e., $\frac{(\tau(i)- \tau(\alpha))}{N} \notin {\Bbb
Z}$. So 
$$
\sum_{g_k} d_{g_k}^2 \omega_{g_k}^{-1} = 
\sum_{1\leq k \leq n} d_g^2 \omega_{g_1}^{-1} z^{k-1} =0
$$. \par
We have 
$$
\align
\sum_{e\in V_0} d_e^2 \omega_e^{-1}
& = \frac{1}{N} \sum_{g\in W_0} d_g^2 \omega_g^{-1} =
\frac{1}{N} \sum_{g\in W} d_g^2 \omega_g^{-1} \\
& = \frac{1}{N} \sum_{k} d_k^2 \omega_k^{-1}
\sum_{\delta} d_\delta^2 \omega_\delta \\
&= \frac{1}{N} \exp(\frac{-6\pi i (C_G-C_H)}{24}) \frac{1}{S_{11}} 
 \frac{1}{\dot{S_{11}}}
\endalign
$$, 
where in the last $=$ we have used (3),(12) and comments after (12) in \S2.1,
and  $C_G$ and $C_H$ are central charges given by (1.4.2) of
[KW], i.e.,
$$
C_G = \frac{(N^2-1)m'}{m'+N} + \frac{(N^2-1)m''}{m''+N},
$$ 
$$
C_H = \frac{(N^2-1)(m'+m'')}{m'+m''+N}.
$$  (2) now follows by a simple calculation.
\hfill Q.E.D \par
\proclaim{Corollary 2.4}
The irreducible covariant sectors of the diagonal coset
of type $A$ corresponding to the basis of $V_0$, with its braiding
and $S,T$ matrices as defined in \S2.1, is
a unitary Modular Tensor Category (cf. [Tu]).
\endproclaim
{\it Proof:}
By the definition of  unitary Modular Tensor Category as on P. 74 
and P. 113 of [Tu], it is enough to show that the  $Y$ matrix 
is invertible, which follows from (1) of Th.2.3. 
\hfill Q.E.D\par
When the $\Bbb {Z}$ action on $exp$ is free, i.e., for 
any $(i,\alpha)\in exp$,   
$\sigma^s (i,\alpha) = (i,\alpha)$ iff $N|s$, it is easy to see  
using (2) of Th.2.3  that  $S_{AB} = N S_{ij} \overline {\dot{S_{\alpha
\beta}}}$ where $A:= P(i\otimes \alpha), B:= P(j\otimes \beta)$ and
both $A,B$ are irreducible. \par
Let us calculate $S$-matrices in the case $N$ is a prime and there exists
a (necessarily unique)
fixed point $F:=P(i_0\otimes \alpha_0) \in V_0$ under the action
of
$\sigma$. By lemma 2.1, 
$F=\sum_{1\leq i \leq N} F_i$ with $F_i$ irreducible and $d_{F_i} = 
\frac{1}{N} d_F$.
Recall $S_{ab} := |\tilde \sigma |^{-1} Y_{ab}$,
where
$$
|\tilde \sigma |^2 
= \frac{1}{N^2} \frac{1}{S_{11}^2} \frac{1}{\dot{{S_{11}^2}}}
$$, which follows from (2) of Th.2.3.  \par
So $S_{ab} = N S_{11}\dot {{S_{11}}} Y_{ab}$, and by (3) of Lemma 2.2
$S_{ab}$ for $a=P(i\otimes \alpha)\neq F_i$ is determined as follows:
$$
S_{a b} = N S_{ij} \overline {\dot {{S_{\alpha \beta}}}} \tag 2
$$ 
if $ b:=P(j\otimes \beta) \neq F_k$, and 

$$
S_{a F_k} =  S_{ii_0} \overline {\dot {{S_{\alpha \alpha_0}}}} \tag 3
$$. \par
It remains to determine $S_{F_iF_j}$.  Note 
$$
STS = T^{-1} S T^{-1}
$$, and  $T_{ab} = \delta_{ab} \dot \omega_a$, where
$\dot \omega_a = \omega_a \omega$, and $\omega$ is 
determined by (2) of Th.2.3 and (3) of \S2.1 as:
$$
\omega:= \exp( \frac{-2\pi i}{24} (C_G-C_H))
$$.  We  have
$$
\align
\dot {\omega_{F_i}}^{-1} S_{F_iF_k} \dot {\omega_{F_k}}^{-1} 
& = \sum_a S_{F_iA} \dot {\omega_{A}} S_{AF_k} \\
&= \sum_l S_{F_iF_l} \dot {\omega_{F_l}} S_{aF_k} 
+ \sum_{A\neq F}  S_{F_iA} \dot {\omega_{A}} S_{AF_k}
\\
&= \dot {\omega_{F}} \sum_l S_{F_iF_l} 
S_{F_lF_k} + \sum_{A\neq F}  S_{F_iA} \dot {\omega_{A}} S_{AF_k}
\\
&= \dot {\omega_{F}} (\delta_{F_i \bar F_k} -
\sum_{B\neq F}  S_{F_i B} S_{ B F_k})                         
 + \sum_{A\neq F}  S_{F_iA} \dot {\omega_{A}} S_{AF_k}
\\
&= \dot {\omega_{F}} \delta_{F_i \bar F_k} 
 + \sum_{A\neq F}  S_{F_iA} (\dot {\omega_{A}}- \dot
{\omega_{F}}) S_{AF_k}\\
&= \dot {\omega_{F}} \delta_{F_i \bar F_k}
 + \frac{1}{N^2} \sum_{A\neq F}  S_{FA} (\dot {\omega_{A}}- \dot
{\omega_{F}}) S_{AF}\\
&= \dot {\omega_{F}} \delta_{F_i \bar F_k}
 + \frac{1}{N^2} \sum_{A}  S_{FA} (\dot {\omega_{A}}- \dot
{\omega_{F}}) S_{AF}
\endalign
$$, where we used $S^2$ is equal to conjugate matrix in the fourth $=$,
and (3) of Lemma 2.2 in the sixth $=$.
So we have
$$
S_{F_iF_k} =\dot {\omega_{F}}^3 + \dot {\omega_{F}}^2 
\frac{1}{N^2} \sum_{A\in V_0}  S_{FA} (\dot {\omega_{A}}- \dot
{\omega_{F}}) S_{AF}
$$.  But
$$
\align
\sum_{A\in V_0}  S_{FA} (\dot {\omega_{A}}- \dot
{\omega_{F}}) S_{AF} 
&= 
  \sum_{A\in V_0, 1\leq i,j\leq N}  S_{F_iA} (\dot
{\omega_{A}}-
\dot
{\omega_{F}}) S_{AF_j} \\
&=   \sum_{A\in V_0, 1\leq i,j\leq N}  S_{F_iA} \dot
{\omega_{A}} S_{AF_j} -  \sum_{A\in V_0, 1\leq i,j\leq N} S_{F_iA}
\dot
{\omega_{F}} S_{AF_j} \\
&=  \sum_{ 1\leq i,j\leq N}  \dot{\omega_F}^{-2} S_{F_iF_j}-  
\sum_{A\in V_0,
1\leq i,j\leq N} S_{F_iA}
\dot
{\omega_{F}} S_{AF_j} \\
&=      \dot{\omega_F}^{-2} S_{FF}- \dot {\omega_{F}} \sum_{
1\leq i,j\leq N} \delta_{F_i \bar F_j}\\
&= \dot{\omega_F}^{-2} S_{FF}- \dot {\omega_{F}} N
\endalign
$$, where $S_{FF} = |\tilde \sigma|^{-1} Y_{FF}$, and in the fourth $=$
we used (4) of lemma 2.2, and 

in the last step we used $\bar F =F$ and 
$$
\sum_{
1\leq i,j\leq N} \delta_{F_i \bar F_j} = \langle F, \bar F \rangle
= \langle F,  F \rangle =N
$$.  So we have:
$$
S_{F_i F_k} = \dot {\omega_{F}}^3 \delta_{F_i\bar F_k} + \frac{1}{N^2}
(S_{FF} - \dot {\omega_{F}}^3 N)
$$. \par
Let us  show 
$$
\dot {\omega_{F}}^3 = {\omega_{F}}^3 \omega^3 = 1
$$. Recall $\omega_F = \exp (2\pi i \Delta_F)$, and since $F$ is the 
unique fixed point, by a simple calculation using (10) of \S2.1
we get 
$$
\Delta_F = \frac{N^2-1}{24} - \frac{N(N^2-1)}{24} ( \frac{1}{m'+N} +
\frac{1}{m'+N} - \frac{1}{m'+m''+N}),
$$
and it follows that  
$$
\omega = \exp (-2\pi i \frac{1}{24} (C_G-C_H)) = \omega_F^{-1}
$$ by (2) of Th. 2.3,  so $\dot {\omega_{F}} = 1$. \par 
Therefore we have
$$
S_{F_i F_k} =  \delta_{F_i\bar F_k} + \frac{1}{N^2}
(S_{FF} -  N)
$$. Since $\bar F=F$, $S_{FA}  $ is real for all $A$, so  
$S_{F_k a}    $ is real for any irreducible $a$, and we must have
$\bar F_k = F_k$ 
since $S$ matrix is invertible. So:

$$
S_{F_i F_k} =  \delta_{ik} + \frac{1}{N^2}
(S_{FF} -  N) \tag 4
$$. 
The formula (2),(3) and (4) above  agree with formula (4.40) of [FSS1] 
(Note
our $S_{FF} = N S_{i_0 i_0} \overline {\dot{S_{\alpha_0 \alpha_0}}}$, where
$F=P(i_0\otimes \alpha_0)$). However one should notice that
our definition of $S$ matrices are very different from those of 
[FSS1]. \par
By formula (7) in \S2.1, we can also determine all the fusion coefficients
in this case using (2),(3) and (4). By (1)
of lemma 2.2, it is enough to determine $N_{F_iF_j}^{F_k}$, and the
formula is:
$$
N_{F_iF_j}^{F_k} = \frac{1}{N^3} \langle FF,F\rangle +
\frac{S_{FF}}{S_{1F}} ( \frac{ \delta_{jk} +  \delta_{ji} +
\delta_{ik}}{N} - \frac{3}{N^2}) + \frac{1}{S_{1F}} (
N \delta_{i,j,k} + \frac{2}{N} - \delta_{jk} -  \delta_{ji} -   
\delta_{ik})
$$, where $\langle FF,F\rangle$ can be calculated by using 
the knowledge of the fusion coefficients of $G$ and $H$ by Th.4.3 of
[X4].  \par
In [FSS1] and [FSS2], certain formula about $S$ matrices 
were derived from other considerations  in the case when $N$ is 
not prime and other types of simple simply connected Lie groups,
and it will be interesting to extend our calculations
above to these cases  and to see if the results agree with 
[FSS1] and [FSS2].\par  
By Cor.2.4, one may calculate 3-manifold invariants, denoted by
$\tau_{G/H}$,  using $S$ matrices
obtained above
as in [Tu]. Our calculations on some lens space suggest that
$$
\tau_{G/H} (M) = \frac{1}{N^{c(M)}} \tau_{G}({M}) \overline \tau_{H}({M})
$$, where $c(M) \in {\Bbb Z}$ depends on the three manifold $M$,
and it will be interesting to see if this is true in 
general. \par
\heading {3. Miscellaneous Results} \endheading
\subheading{3.1 KWH and KWC} 
Let $H\subset G_k$ be as in  the introduction. 
Through
out this section, we will assume the following:
$H$ and $G$ verifies similar statements as Cor. 1 in \S34 of [W2] 
(cf: comments after (12) in \S2.1) , and
$H\subset G_k$ is cofinite as defined on Page 18 of [X4]. \par
Note the assumption is satisfied by many examples (cf. Cor. 4.2 of [X4])
and is expected to be true in general.\par
We also assume that  $H\subset G_k$ is not conformal, so the coset 
theory is non-trival. \par
We will use the notations of \S4.2 of [X4] and ideas of [X2]. 
We  denote the set of irreducible sectors of $\sigma_i a_{1\otimes 
\lambda}$  by $V$.  Notice $\sigma_i \in V$,
and these are referred to as "special nodes" in \S3.4 of [X1].  Let:
$$           
 a_{1\otimes\lambda}  a
= \sum_{b\in V} V^\lambda_{ab} b,
$$ where $V^\lambda_{ab}$ are nonnegative
integers.  Denote by $V^\lambda$ the matrix such that $(V^\lambda)_a^b = 
V^\lambda_{ab}$.
Define matrix $N_c$ by $N_{ca}^b = \langle ca,b\rangle$ for
 $a,b,c \in V$.
Then $V^\lambda = \sum_c V^\lambda_{1c} N_c$.
Since $[a_{1\otimes \bar\lambda}] = [\bar a_{1\otimes \lambda}], [\sigma_j
a_{1\otimes \lambda}]
= [a_{1\otimes \lambda} \sigma_j], \ V^\lambda, \ N_{\sigma_j}$ are
commuting normal matrices, so they can be simultaneously diagonalized.
Recall the irreducible representations of the ring $Gr(C_k)$ generated
by $\lambda's$ are given by
$$           
\lambda \rightarrow \frac{S_{\lambda \mu}}{S_{1\mu}}.
$$           
Assume       
$$           
V^\lambda_{ab} = \sum_{i,\mu,s\in(\text{\rm Exp})}
\frac{S_{\lambda \mu}}{S_{1\mu}}  \psi_a^{(i,\mu,s)}
\psi_b^{(i,\mu,s)^*}
$$,  where $ \psi_a^{(i,\mu,s)}$ are normalized
orthogonal   
eigenvectors of $V^\lambda$ (resp. $N_{\sigma_i}$) with eigenvalue
$\frac{S_{\lambda \mu}}{S_{1\mu}}$ (resp.
 $\frac {S_{i j}}{S_{1 j}}$) .
$(Exp)$ is a set of $i,\mu,s$'s
 and $s$ is an
index indicating the multiplicity of  $i,\mu$. We denote by $Exp$
the set of $(i,\mu)$ such that $(i,\mu,s) \in {(Exp)}$ for
some $s$.    
Recall if a representation is denoted by $1$, it will always be the
vacuum representation.  The Perron-Frobenius eigenvector $\psi^{(1,1)}$
is given by $\sum_a d_a a$, up to a positive constant. Note 
all the entries of $\psi^{(1,1)}$ are positive.
\proclaim{Proposition 3.1}
$(i,\alpha)\in {Exp}$ if and only if $b(i,\alpha)>0$.  
\endproclaim 
\noindent    
{\it Proof:} 
Recall (1) of \S1:
$$           
b(i,\alpha) = \sum_{(j,\beta) } S_{ij} \overline{\dot S_{\alpha \beta}}
\langle (i,\alpha), (1,1) \rangle.
$$           
By the proof of (2) of Prop. 4.2 of [X4] and (2) of Cor. 3.5  of 
[X1]  we     
have         
$$           
\langle (i,\alpha), (1,1) \rangle = \langle \sigma_i a_{1\otimes
\bar{\alpha}}, 1 \rangle = \langle \sigma_{\bar i} a_{1\otimes
{\alpha}}, 1 \rangle
$$, so:      
$$           
\align       
b(i,\alpha) & = \sum_{(j,\beta) } S_{ij} \overline{\dot S_{\alpha \beta}}
\langle \sigma_{\bar j} a_{1\otimes
{\beta}}, 1 \rangle = \sum_{(j,\beta) } S_{ij} \overline{\dot S_{\alpha
\beta}} (N_{\bar{j}} V^{\beta})_{11} \\
& = \sum_{ (k, \delta, s) \in ({Exp})} 
S_{ij} \overline{\dot S_{\alpha \beta}} \frac{\dot S_{\beta \delta}}{ 
\dot S_{1
\delta}}     
\frac{\overline{S_{jk}}}{ S_{1 k}} |\psi_1^{(k, \delta, s)}|^2 \\
&= \sum_{s} 
\frac{1}{\dot S_{1 \alpha}} \frac{1}{S_{1i}} 
 |\psi_1^{(i,\alpha,s)}|^2 
\endalign    
$$.          
Note the equality above is similar to (1) on Page 12 of [X2], and
the rest of the proof is the same as the proof on Page 12 of 
[X2].        
\hfill Q.E.D. \par 
Note if $b(i,\alpha) >0$, then $(i,\alpha) \in exp$, so by Prop. 3.1
${Exp} \subset exp$, and KWC is equivalent to the statement that
${Exp} = exp$. \par
\proclaim{Proposition 3.2}
If $b(i,1)>0$, then $\frac{b(i,1)}{b(1,1)} = \frac{S_{i1}}{S_{11}}$.
\endproclaim 
{\it Proof:} 
If  $b(i,1)>0$, by Prop.3.1, $(i,1) \in {Exp}$. Suppose
$(j,\beta) \succ (1,1)$. Then
$$           
\langle (j,\beta) , (1,1) \rangle = \langle \sigma_j a_{1\otimes \bar
\beta}, 1 \rangle = \langle \sigma_j,  a_{1\otimes \beta} \rangle >0
$$, since $\sigma_j$ is irreducible, it follows that
$$           
a_{1\otimes \beta} \succ \sigma_j
$$, and if  $b= a_{1\otimes \beta} - \sigma_j$, then $V^b:= V^\beta - N_j
$ is a normal
matrix with non-negative entries, with a Perron-Frobenius eigenvalue
$ \frac{\dot S_{\beta 1}}{\dot S_{11}} - \frac{S_{j1}}{S_{11}}$. It
follows that 
$$           
\align       
\frac{\dot S_{\beta 1}}{\dot S_{11}} - \frac{S_{j1}}{S_{11}} & \geq |
\langle V^b \psi^{(i,1,s)}, \psi^{(i,1,s)} \rangle |= 
|\frac{\dot S_{\beta 1}}{\dot S_{11}} - \frac{S_{ji}}{S_{1i}}| \\
& \geq \frac{\dot S_{\beta 1}}{\dot S_{11}} - |\frac{S_{ji}}{S_{1i}}| \\
& \geq \frac{\dot S_{\beta 1}}{\dot S_{11}} - \frac{S_{j1}}{S_{11}}
\endalign    
$$, and so all the $\geq$ are $=$, which  happens only if
$$           
\frac{S_{ji}}{S_{1i}} = \frac{S_{j1}}{S_{11}}
$$. Prop. 3.1 now follows from the definitions.
\hfill Q.E.D. 
             
\proclaim{Theorem 3.3}
If $H\subset G_k$  satisfies (2) in the introduction,   then KWH is
equivalent to KWC.
\endproclaim 
{\it Proof:} 
We just have to show that KWC implies KWH. \par
By (2) of the introdcution we have
if  $\langle (1,\delta), (1,1) \rangle >0$, then $\delta =1$. For
any $\alpha, \beta$, we have (cf. Prop.4.2 of [X4]):
$$           
\align       
\langle a_{1\otimes \alpha}, a_{1\otimes \beta} \rangle & =
\langle a_{1\otimes \alpha} a_{1\otimes \bar \beta}, 1 \rangle
=  N_{\alpha \bar\beta}^\delta \langle  a_{1\otimes \delta}, 1 \rangle \\
&=  N_{\alpha \bar\beta}^\delta \langle (1,\delta), (1,1) \rangle
\\           
&=  N_{\alpha \bar\beta}^1 \\
&= \langle \alpha, \beta \rangle
\endalign    
$$. In particular $a_{1\otimes \beta}, \forall \beta$ is irreducible.
\par         
Suppose $\langle (j,\beta), (1,1) \rangle >0$, then 
$\langle \sigma_j a_{1\otimes \bar \beta}, 1 \rangle >0$, and so
$\langle \sigma_j , a_{1\otimes  \beta} \rangle >0$. Since both 
$\sigma_j$ and $a_{1\otimes  \beta}$ are irreducible, it follows that
$$           
\sigma_j = a_{1\otimes  \beta}
$$.\par      
Now suppose $(i,\alpha)\in exp$. By KWC, $b(i,\alpha) >0$, 
so  by Prop. 3.1, $(i,\alpha) \in {Exp}$, and 
from 
$\sigma_j = a_{1\otimes  \beta}$ we 
must have:   
$$           
\frac{S_{ji}}{S_{1i}} = \frac{\dot S_{\beta \alpha}}{\dot S_{1 \alpha}} 
$$, and therefore
$$           
S_{ji} \overline{\dot S_{\beta \alpha}} = S_{1i} \frac{|\dot S_{\beta
\alpha}|^2}{\dot S_{1
\alpha}} \geq 0
$$, which proves KWH.
\hfill Q.E.D. \par
Let us give an example which does not satisfy the assumption 
of our theorem, and verifies KWC but not KWH.  This is the
coset $SU(2)_8 \subset SU(3)_2$
discussed on Page 32 of [X4] (also cf. [DJ]) and we will use the notations
there.       
The Vacuum representation space $H$ of $LSU(3)$ at level $2$ decomposes
as:          
$$           
H= (00,0) \otimes 0 + (00,4)\otimes 4 + (00,8) \otimes 8
$$, and since 
$$           
(00,0) = (00,8)
$$           
as representations of the coset, our assumption is not satisfied.
By the ring structure given on Page 32 of [X4] we checked C2 is 
satisfied, which implies that KWC is true. However, since 
$$           
(11,4) =(00,0)
$$ and $(10,4)\in exp$, KWH implies that
$$           
S_{11,10} \dot{S_{4,4}} \geq 0
$$.  But a  direct calculation shows 
$$           
\frac{S_{11,10}}{S_{00,10}} = \frac{1-\sqrt 5}{2}
$$ and       
$$           
\dot{S_{4,4}}= \sqrt {\frac{1}{5}}
$$, which shows that
$$           
S_{11,10} \dot{S_{4,4}} <0
$$ since $S_{00,10} >0$.  In fact this was discovered when we verified
C2 in this example. \par   
             
Note that all diagonal inclusions of type $A$ satifies the 
assumption of Th.3.3  by 2.7.12 of
[KW].  To give a slightly different example, let us consider
the following inclusions
$$           
SU(2)_{11k} \subset SU(2)_{8k} \times SU(2)_{3k} \subset SU(3)_{2k} \times
SU(2)_{3k}   
\subset SU(6)_k
$$ with $k\in {\Bbb N}$, where the first inclusion is diagonal, the 
second inclusion comes from the conformal inclusion $SU(2)_4 \subset SU(3)_1$,
and the third inclusion comes from the conformal inclusion $SU(3)_2
\times SU(2)_3 \subset SU(6)$ (cf. [X3]). 

By (2) of Prop.3.1 and (2) of Cor.3.1 of [X4],
the inclusion $SU(2)_{11k} \subset SU(6)_k$ is cofinite and verifies the
assumptions at the beginning of this section. Note if $(1,\alpha) \succ
(1,1)$, then $\Delta_\alpha \in {\Bbb Z}$. Here $ 1\leq \alpha \in 
{\Bbb Z}  \leq 11k+1$ and $\Delta_\alpha = \frac{\alpha (\alpha
-1)}{11k+2}$. So if $11k+2$ is prime, then $(1,\alpha) \succ
(1,1)$ iff $\alpha=1$, and it follows that the conditions of Th.2.3
are satisfied. So the conclusions of Th.2.3 hold for the 
inclusion    
$$SU(2)_{11k} \subset  SU(6)_k$$ if $11k+2$ is prime, and by 
Dirichlet's theorem, there are infinitely many such $k$'s. \par
\subheading{3.2. A property of C2}
\proclaim{Proposition 3.4}
If $H_1\subset H_2 $, $H_2 \subset G_k$ verifies C2 or is a
conformal inclusion, assume $H_1\subset G_k$ is not a conformal inclusion
 to avoid trivality, then
$H_1\subset G_k$ verifies C2.
\endproclaim 
{\it Proof:} 
For simplicity we will use $x,y,z$ to denote the
irreducible representations of $LH_1$, $LH_2$ and $LG$ respectively,
and $A, B ,C$ to denote the Vacumm sector of cosets $H_1\subset H_2,
H_2\subset G, H_1\subset G$ respectively. 
Note we have natural inclusions $A\otimes B\subset C$.
From the decompositions:
$$           
\pi_z \simeq \sum_y \pi_{(z,y)} \otimes \pi_y \simeq \sum_{y, x}
\pi_{(z,y)}  
\otimes \pi_{(y,x)} \otimes \pi_{x} 
\simeq  \sum_{x} \pi_{(z,x)} \otimes \pi_x  
$$           
we conclude that 
$$           
\pi_{(z,x)} \simeq  \sum_y \pi_{(z,y)} \otimes \pi_{(y,x)}
$$, which is understood as the decomposition of representation
$\pi_{(z,x)}$ of $C$ when restricted to $A\otimes B \subset C$.
By local equivalence, the statistical dimension of 
$$           
\pi_{(z,x)}(A(I)\otimes B(I)) \subset \pi_{(z,x)}(C(I))
$$ is the same as  that of 
$$           
\pi_{(1,1)}(A(I)\otimes B(I)) \subset \pi_{(1,1)}(C(I))
$$ ( $I$ is a proper interval of the circle),  which is easily seen 
(by using Haag duality) to be 
$$           
\sum_y d_{(1,y)} d_{(y,1)}
$$. 
Here when  $H_2 \subset G$ (resp. $H_1 \subset H_2$) is a conformal inclusion,
$d_{(z,y)}$ (resp. $d_{(y,x)}$)  is defined to be the multiplicity
of irreducible representation $y$ (resp. $x$)  which appears in $z$
(resp. $y$) when restricting to $LH_2$ (resp. $LH_1$). \par
So the  statistical dimension of the inclusion:
$\pi_{(z,x)}(A(I)\otimes B(I)) \subset \pi_{(z,x)}(C(I)) 
\subset \pi_{(z,x)}(C(I^c))' \subset \pi_{(z,x)}(A(I^c)\otimes B(I^c))'
$ ($I^c$ is the complement of $I$ is $S^1$) is
$$           
d_{(z,x)} \sum_y d_{(1,y)} d_{(y,1)}
$$. On the other hand, by [L6], the statistical dimension of the
above inclusion is also given by:
$$           
\sum_y d_{(z,y)} d_{(y,x)}
$$, and so:  
$$           
d_{(z,x)} \sum_y d_{(1,y)} d_{(y,1)} = \sum_y d_{(z,y)} d_{(y,x)}
$$, therefore 
$$           
d_{(z,x)} =\frac{\sum_y d_{(z,y)} d_{(y,x)}}{\sum_y d_{(1,y)} d_{(y,1)}}
$$. The proof now follows from the assumptions and
$$           
b(z,x) = \sum_y b(z,y)b(y,x)
$$ which follows from Th. B of [KW].
\hfill Q.E.D. \par
Let us give one application of the above proposition. Consider the 
superconformal coset models (cf. [Gep], [LVW] or [NS]):
$$           
G(m,n,k):= \frac{SU(m+n)_k \times SO(2mn)_1}{SU(m)_{n+k} \times
SU(n)_{m+k} \times U(1)_{mn(m+n)(m+n+k)}}
$$.          
In our setting, when $2mn>2$, the inclusion is given by $H\subset G$
with         
$H=SU(m)_{n+k} \times
SU(n)_{m+k} \times U(1)_{mn(m+n)(m+n+k)}
$ and $G= SU(m+n)_k \times Spin(2mn)_1$. The inclusion 
$H\subset G$ is constructed by the composition of two inclusions:
$$   
\align        
H  & \subset SU(m)_{n} \times SU(m)_{k}  \times SU(n)_{m} \times SU(n)_{k} \\ 
& \times U(1)_{mn(m+n)(m+n)} \times U(1)_{mn(m+n)(k)}  \tag 1
\endalign
$$ 
and
$$  
\align
(SU(m)_{n} & \times SU(n)_{m}  
\times U(1)_{mn(m+n)(m+n)})   \times
(SU(m)_{k} \times SU(n)_{k} \\
 &  \times U(1)_{mn(m+n)(k)})   
  \subset  Spin(2mn)_1 \times SU(m+n)_k \tag 2
\endalign
$$. Here we need to explain the inclusion 
$$           
SU(m)_{n} \times SU(n)_{m} \times U(1)_{mn(m+n)(m+n)} \subset 
Spin(2mn)_1  
$$. This inclusion comes from the action of $SU(m) \times SU(n)
\times U(1)$ on the tangent space ($\simeq {\Bbb C}^m \otimes {\Bbb C}^n$)
of Grassmanian 
$$           
\frac{SU(m+n)}{SU(m) \times
SU(n) \times U(1)}
$$, which  gives an inclusion 
$$           
SU(m)_{n} \times SU(m)_{n} \times U(1)_{mn(m+n)(m+n)} \subset
SO(2mn)      
$$ .  But it is easy to check that the above inclusion maps all the
elements of  
the fundamental group of $U(1)$  to trival element in the fundamental
group of $SO(2mn)$ 
(note $SU(m), SU(n)$ are simply connected), 
so the inclusion above lifts to an inclusion into 
$Spin(2mn)$. 
The inclusion 
$$
SU(m)_k \times SU(n)_k \times U(1)_{mn(m+n)(k)} \subset SU(m+n)_k
$$ comes from the conformal inclusion 
$$
SU(m)_1 \times SU(n)_1 \times U(1)_{mn(m+n)} \subset SU(m+n)_1.
$$

\par 
             
The  inclusion in (1) is diagonal, and the $SU$ part of the
inclusion verifies C2 by (4) of Th. 4.3 in [X4]. For the $U(1)$ part,
we consider the following inclusions:
$$           
U(1)_{2a} \times U(1)_{2b} \subset SU(2)_{a} \times SU(2)_{b},
$$ 
and 
$$
U(1)_{2a+2b} \subset SU(2)_{a+b} \subset SU(2)_{a} \times SU(2)_{b},
$$
with $a:= \frac{1}{2} mn (m+n)^2, b:= \frac{1}{2} mn (m+n)k$.
It follows from (3) of Cor.3.1 of [X4] and the proof of (1) of Prop.3.1
that         
$$           
 U(1)_{2a+2b} \subset    U(1)_{2a} \times U(1)_{2b}
$$ is cofinite, and a much simpler argument (since all the sectors
involved are automorphisms) as in the proof of (4) of
Th. 4.3 in [X4]  verifes C2 in this case. \par
It follows from Prop.3.1 and Th.4.2 of [X4], [W2] and 
[JB]  that $G(m,n,k)$ coset
verifies Conj.1 of [X4], and so is indeed a "rational" conformal 
field theory.
\par         
By Proposition 3.4, we see that when k=1 and $mn>1$, the above coset
verifies C2.\par 
The fixed point resolution problems for 
$G(m,n,k)$ are discussed in [Gep], [LVW]
(also        
cf. [NS]). It will be interesting to work out this problem along the lines 
of \S2.      
\heading References \endheading
\roster
\item"[DJ]" D. Dunbar and K. Joshi,
{\it Characters for Coset conformal field theories and Maverick
examples}, Inter. J. Mod. Phys. A, Vol.8, No. 23 (1993), 4103-4121.
\item"[Dyn1]" E. B. Dynkin, Amer. Math. Soc. Transl. (2) , vol.6, 
245, (1957).
\item"[Dyn2]" E. B. Dynkin, Amer. Math. Soc. Transl. (2) , vol.6, 
111, (1957).
\item"{[FRS]}" K.Fredenhagen, K.-H.Rehren and B.Schroer
,\par
{\it Superselection sectors with braid group statistics and 
exchange algebras. II}, Rev. Math. Phys. Special issue (1992), 113-157.
\item"{[FSS1]}" J. Fuchs, B. Schellekens and C. Schweigert,
{\it The resolution of field identification fixed points in diagonal coset
theories,} Nucl. Phys. B 461 (1996) 371, hep-th/9509105.
\item"{[FSS2]}" J. Fuchs, B. Schellekens and C. Schweigert,
{\it From Dynkin diagram symmetries to fixed point structures,}
Commun. Math. Phys. 180 (1996) 39
\item"{[Gep]}" D. Gepner, Phys. Lett. B222 (1989) 207 
\item"{[GL]}"  D.Guido and R.Longo, {\it  An Algebraic Spin and
Statistics Theorem},  \par
Comm.Math.Phys., {\bf 181}, 11-35 (1996)  
\item"{[GL2]}"  D.Guido and R.Longo, {\it  Relativistic invariance and
charge
conjugation in quantum field theory},  \par
Comm.Math.Phys., {\bf 148}, 521-551 (1992)
\item"{[JB]}" Jens Bockenhauer, {\it an algebraic formulation of level 1
WZW models,}, Hep-th 9507047, Rev.Math.Phys. {\bf 8}, 925-948 (1996).
\item"{[J]}" V. Jones, {\it Fusion en al\'gebres de Von Neumann et groupes
de lacets (d'apr\'es A. Wassermann),} Seminarie Bourbaki, 800, 1-20,1995.
\item"{[KW]}"  V. G. Kac and M. Wakimoto, {\it Modular and conformal
invariance constraints in representation theory of affine algebras},  
Advances in Math., {\bf 70}, 156-234 (1988).
\item"{[Kac]}"  V. G. Kac, {\it Infinite dimensional algebras}, 3rd
Edition,
Cambridge University Press, 1990. 
\item"{[L1]}"  R. Longo, Proceedings of International Congress of
Mathematicians, 1281-1291 (1994).
\item"{[L2]}"  R. Longo, {\it Duality for Hopf algebras and for
subfactors}, 
I, Comm. Math. Phys., {\bf 159}, 133-150 (1994).
\item"{[L3]}"  R. Longo, {\it Index of subfactors and statistics of
quantum fields}, I, Comm. Math. Phys., {\bf 126}, 217-247 (1989.
\item"{[L4]}"  R. Longo, {\it Index of subfactors and statistics of
quantum fields}, II, Comm. Math. Phys., {\bf 130}, 285-309 (1990).
\item"{[L5]}"  R. Longo, {\it An analog of the Kac-Wakimoto formula and
black hole conditional entropy,} gr-qc 9605073, to appear in
Comm.Math.Phys.
\item"{[L6]}"  R. Longo, {\it Minimal index and braided subfactors,
} J.Funct.Analysis {\bf 109} (1992), 98-112.
\item"{[LR]}"  R. Longo and K.-H. Rehren, {\it Nets of subfactors},
Rev. Math. Phys., {\bf 7}, 567-597 (1995).  
\item"{[LVW]}" W. Lerche, C. Vafa and N. P. Warner, Nucl. Phys. B324
(1989) 427.
\item"[NS]" S. Naculich and H. Schnitzer, {\it Superconformal 
coset equivalence from level-rank duality}, hep-th/9705149.
\item"{[PP]}" M.Pimsner and S.Popa,
{\it Entropy and index for subfactors}, \par
Ann. \'{E}c.Norm.Sup. {\bf 19},
57-106 (1986). 
\item"[PS]" A. Pressly and G. Segal, {\it Loop Groups,} O.U.P. 1986.
\item"[Reh]" Karl-Henning Rehren, {\it Braid group statistics and their
superselection rules} In : The algebraic theory of superselection 
sectors. World Scientific 1990 
\item"[SY]" A. N. Schellekens and S. Yankielowicz, Nucl. Phts. B 324,
67, (1990).  
\item"[Tu]" V. G. Turaev, {\it Quantum invariants of knots and
3-manifolds,} Walter de Gruyter, Berlin, New York 1994.
\item"{[W1]}"  A. Wassermann, Proceedings of International Congress of
Mathematicians, 966-979 (1994).
\item"{[W2]}"  A. Wassermann, {\it Operator algebras and Conformal
field theories III},  Invent. Math. Vol. 133, 467-539 (1998)
\item"[X1]" F.Xu, {\it   New braided endomorphisms from conformal
inclusions, } \par
Comm.Math.Phys. 192 (1998) 349-403.
\item"[X2]" F.Xu, {\it Applications of Braided endomorphisms from
Conformal inclusions,} 
Inter. Math. Res. Notice., No.1, 5-23 (1998), see also q-alg/9708013,
and Erratum, Inter. Math. Res. Notice., No.8, (1998) 
\item"[X3]" F.Xu, {\it Jones-Wassermann subfactors for 
Disconnected Intervals}, Preprint 97, see also
 q-alg/9704003.
\item"[X4]" F.Xu, {\it Algebraic coset conformal field theories},
preprint 98, q-alg/9810053.
\endroster   
\enddocument